\documentclass[12pt,a4paper,reqno]{amsart}
\usepackage{amsmath}
\usepackage{amsthm}
\usepackage{amsfonts}
\usepackage{amssymb}
\usepackage{setspace}
\numberwithin{equation}{section}

\theoremstyle{plain}

\newtheorem{theorem}[subsection]{Theorem}
\newtheorem{proposition}[subsection]{Proposition}
\newtheorem{lemma}[subsection]{Lemma}

\newtheorem{conjecture}[subsection]{Conjecture}

\newtheorem{definition}[subsection]{Definition}

\newtheorem*{a1}{Lemma \ref{prog-lem}}

\theoremstyle{definition}

\newtheorem{remark}[subsection]{Remark}

\renewcommand{\leq}{\leqslant}
\renewcommand{\geq}{\geqslant}

\newsavebox{\proofbox}
\savebox{\proofbox}{\begin{picture}(7,7)%
 \put(0,0){\framebox(7,7){}}\end{picture}}

\newcommand{\md}[1]{\ensuremath{(\mbox{mod}\, #1)}}

\newcommand{\mdlem}[1]{\ensuremath{(\mbox{\textup{mod}}\, #1)}}

\def\E{\mathbb{E}}
\def\Z{\mathbb{Z}}
\def\R{\mathbb{R}}

\def\C{\mathbb{C}}
\def\N{\mathbb{N}}
\def\P{\mathbb{P}}
\def\F{\mathbb{F}}

\newcommand\PFR{\operatorname{PFR}}
\newcommand\PGI{\operatorname{PGI}}
\newcommand\rk{\operatorname{rk}}

\def\eps{\varepsilon}

\DeclareMathOperator{\hcf}{hcf}

\def\proof{\noindent\textit{Proof. }}
\def\endproof{\hfill{\usebox{\proofbox}}\vspace{11pt}}

\begin{document}

\onehalfspace
\title[Equivalence of PFR and PGI(3)]{An equivalence between inverse sumset theorems and inverse conjectures for the $U^3$ norm}

\author{Ben Green}
\address{Centre for Mathematical Sciences\\
Wilberforce Road\\
     Cambridge CB3 0WA\\
     England
}
\email{b.j.green@dpmms.cam.ac.uk}

\author{Terence Tao}
\address{UCLA Department of Mathematics, Los Angeles, CA 90095-1555
}
\email{tao@math.ucla.edu}

\thanks{The first author holds a Leverhulme Prize and is grateful to the Leverhulme Trust for their support. The second author is supported by a grant from the MacArthur Foundation, and by NSF grant DMS-0649473.  This work was completed while the authors were attending the programme on \emph{Ergodic Theory and Harmonic Analysis} at MSRI and they would like to thank the Institute for providing excellent working conditions.}

\begin{abstract}
We establish a correspondence between \emph{inverse sumset theorems} (which can be viewed as classifications of approximate (abelian) groups) and \emph{inverse theorems for the Gowers norms} (which can be viewed as classifications of approximate polynomials).  In particular, we show that the inverse sumset theorems of Fre\u{\i}man type are equivalent to the known inverse results for the Gowers $U^3$ norms, and moreover that the conjectured polynomial strengthening of the former is also equivalent to the polynomial strengthening of the latter.  We establish this equivalence in two model settings, namely that of the finite field vector spaces $\F_2^n$, and of the cyclic groups $\Z/N\Z$.

In both cases the argument involves clarifying the structure of certain types of \emph{approximate homomorphism}.
\end{abstract}
\maketitle

\section{introduction}

\textsc{Approximate groups.} The notion of an approximate group has come to be seen as a central one in additive combinatorics. Let $K \geq 1$ be a parameter (the ``roughness'' parameter), and suppose that $A$ is a finite subset of some ambient abelian group $G = (G,+)$ (such as the integers $\Z$). We say that $A$ is a \emph{$K$-approximate group} if $A$ is symmetric (that is to say $-x \in A$ whenever $x \in A$) and if the sumset $A + A := \{a + a': a,a' \in A\}$ is covered by $K$ translates of $A$.  Thus, for instance, the arithmetic progression $\{-N,\ldots,N\}$ in the integers $\Z$ for any $N \geq 1$ is a $3$-approximate group, while the $1$-approximate group are nothing more than the finite subgroups of $G$.

The basic theory of approximate abelian groups was developed by Ruzsa in several papers \cite{ruzsa-freiman,ruzsa,ruzsa:f2}; see also \cite{tao:approxgroups} for some extensions to non-abelian groups. 

Perhaps the most basic question to ask about an approximate group is that of the extent to which it resembles an actual group. A language for formalising this was introduced by the second author in \cite{tao:solvable}, and in the abelian case it reads as follows.

\begin{definition}[Control]
Let $A$ and $B$ be two sets in some ambient abelian group, and $K \geq 1$. We say that $B$ \emph{$K$-controls} $A$ if $|B| \leq K|A|$ and if there is some set $X$ in the ambient group with $|X| \leq K$ and such that $A \subseteq B + X$.
\end{definition}

Two of the landmark results of additive combinatorics may be stated in this language. The first of these may be found in \cite{ruzsa:f2} and the second in \cite{chang-freiman}, a paper which builds upon \cite{frei} and \cite{ruzsa-freiman}.

\begin{theorem}[Inverse sumset theorem for $\F_2^\infty$]\label{inv1}
Suppose that $A \subseteq \F_2^{\infty}$ is a $K$-approxi-mate group for some $K \geq 2$. Then $A$ is $e^{K^C}$-controlled by a \textup{(}genuine\textup{)} finite subgroup of $\F_2^\infty$. 
\end{theorem}

\begin{theorem}[Inverse sumset theorem for $\Z$]\label{inv2}
Suppose that $A \subseteq \Z$ is a $K$-approximate group for some $K \geq 2$. Then $A$ is $e^{K^C}$-controlled by a symmetric generalized arithmetic progression $P = \{l_1 x_1 + \dots + l_d x_d : l_i \in \Z, |l_i| \leq L_i \hbox{ for all } 1 \leq i \leq d\}$ with dimension $d \leq K^C$. 
\end{theorem}

\begin{remark} In this paper the letter $C$ will always denote a constant, but different instances of the notation may indicate different constants.  The restriction $K \geq 2$ is purely a notational convenience, so that we may write $K^C$ instead of $CK^C$.
\end{remark}

These theorems, the background to them and their proofs are now discussed in many places. See, for example, the book \cite{tao-vu}. Neither result is usually formulated in precisely this fashion, but simple arguments involving the covering lemmas in \cite[Chapter 2]{tao-vu} may be used to deduce the above forms from the standard ones.  The proofs of the above two theorems extend easily to the case of bounded torsion $G$ and torsion-free $G$ respectively.  It is also possible to establish a result valid for all abelian groups at once, and containing the above two results as special cases: see \cite{green-ruzsa} for details.

There seems to be a general feeling that the bounds in these results are not optimal, and the so-called \emph{Polynomial Fre\u{\i}man-Ruzsa conjecture} (PFR) has been proposed as a suggestion for what might be true.

\begin{conjecture}[PFR over $\F_2^{\infty}$]\label{PFR-2}
Suppose that $A \subseteq \F_2^{\infty}$ is a $K$-approximate group. Then $A$ is $K^C$-controlled by a finite subgroup.
\end{conjecture}

\begin{conjecture}[Weak PFR over $\Z$]\label{PFR-1}
Suppose that $A \subseteq \Z$ is a $K$-approximate group. Then $A$ is $e^{K^{o(1)}}$-controlled by a symmetric generalised arithmetic progression $P = \{l_1 x_1 + \dots + l_d x_d : |l_i| \leq L_i\}$  with dimension $d \leq K^{o(1)}$, where $o(1)$ denotes a quantity bounded in magnitude by $c(K)$ for some function $c$ of $K$ that goes to zero as $K \to \infty$.
\end{conjecture}

Conjecture \ref{PFR-2} has been stated in several places, and in the article \cite{green:finfield} unpublished work of Ruzsa was discussed, establishing a number of equivalent forms of it. According to Ruzsa \cite{ruzsa-freiman}, the first person to make a conjecture equivalent to the PFR over $\F_2^{\infty}$ was Katalin M\'arton. Conjecture \ref{PFR-1}, concerning approximate subgroups of $\Z$, does not to our knowledge appear explicitly in the literature, although something close to it was suggested by Gowers \cite{gowers:gafa}. One might very optimistically conjecture that a $K$-approximate subgroup of $\Z$ is $K^C$-controlled by the affine image of the set of lattice points inside a convex body of dimension $O(\log K)$. Such a conjecture might deserve to be called the PFR over $\Z$ (rather than the \emph{weak} PFR), since it is nontrivial even if $K$ is a suitably small power of $|A|$. A number of issues are rather unclear concerning such a formulation, one of them being whether it suffices to consider \emph{boxes} rather than arbitrary convex bodies. This question appears to involve somewhat subtle issues from convex geometry and we will not consider it, or indeed any aspect of the stronger version of the PFR over $\Z$, any further in this paper.\vspace{11pt}

\textsc{approximate polynomials} We turn now to what appears to be a completely unrelated topic. Let $G = (G,+)$ be a finite abelian group, and recall the definition of the Gowers norms. If $f : G \rightarrow \C$ is a function we define
\begin{align*}
\Vert f \Vert_{U^1(G)} &:= (\E_{x,h \in G} f(x) \overline{f(x+h)})^{1/2} \\
\Vert f \Vert_{U^2(G)} &:= (\E_{x,h_1,h_2 \in G} f(x) \overline{f(x+h_1)f(x+h_2)}f(x+h_1+h_2))^{1/4} \\
\Vert f \Vert_{U^3(G)} &:= (\E_{x,h_1,h_2,h_3 \in G} f(x) \overline{f(x+h_1)f(x+h_2)f(x+h_3)}\times \\ 
&\quad \times f(x+h_1+h_2)f(x+h_1+h_3)f(x+h_2+h_3)\overline{f(x+h_1+h_2+h_3)})^{1/8}
\end{align*} 
and so forth, where we use the averaging notation $\E_{x \in A} f(x) := \frac{1}{|A|} \sum_{x \in A} f(x)$. In this paper we shall be working primarily with the $U^3(G)$-norm. It is clear that if $\Vert f \Vert_{\infty} \leq 1$ and $\Vert f \Vert_{U^3(G)} = 1$ then we necessarily have $f(x) = e(\phi(x))$, where $\phi : G \rightarrow \R/\Z$ is a \emph{quadratic polynomial} in the sense that $\Delta_{h_1} \Delta_{h_2} \Delta_{h_3} \phi(x)  = 0$ for all $h_1,h_2,h_3,x \in G$, where $\Delta_h \phi(x) := \phi(x+h)-\phi(x)$.  To justify the terminology, observe that when $G = \Z/N\Z$ with $N$ odd it is an easy matter to check that any quadratic polynomial has the form $\phi(x) = \frac{1}{N} ax^2 + \frac{1}{N} bx + c$ for $a,b \in \Z/N\Z$ and $c \in \R/\Z$, where $\frac{1}{N}: \Z/N\Z \to \R/\Z$ is the usual embedding.

The \emph{inverse problem for the Gowers $U^3$-norm} asks what can be said about functions $f : G \rightarrow \C$ for which $\Vert f \Vert_{\infty} \leq 1$ and $\Vert f \Vert_{U^3(G)} \geq 1/K$. In view of the above discussion, it is reasonable to call such functions $f$ \emph{$K$-approximate quadratics.} 

The analogue of \emph{control} in this setting is \emph{correlation}. We say that a function $f : G \rightarrow \C$ \emph{$\delta$-correlates} with another function $F : G \rightarrow \C$ if the inner
product $\langle f, F\rangle := \E_{x \in G} f(x) \overline{F(x)}$ is at least $\delta$.

In the finite field setting, the following inverse theorem was shown in \cite{samorodnitsky}.

\begin{theorem}[Inverse theorem for $U^3(\F_2^n)$]\label{invu3-f2}
Suppose that $f : \F_2^n \rightarrow \C$ is a $K$-approximate quadratic for some $K \geq 2$. Then $f$ $\exp(-K^C)$-correlates with a phase $(-1)^{\psi}$ for some quadratic polynomial $\psi: \F_2^n \to \F_2$.
\end{theorem}

\emph{Remark.} The phase $\psi(x)$ may be written explicitly, relative to a basis, as $\psi(x) := x \cdot Mx + b \cdot x + c$, where $M: \F_2^n \to 
\F_2^n$ is a linear transformation, $b \in \F_2^n$, $c \in \F_2$, and $b \cdot x$ is the usual dot product in $\F_2$.\vspace{11pt}

In $\Z/N\Z$ there is an analogous result, which we recall in Theorem \ref{invu3} below. To state it we recall some of the terminology from \cite{green-tao-u3inverse} concerning \emph{nilsequences}.

\begin{definition}[Nilsequences] A \emph{$2$-step nilmanifold} is a homogeneous space $G/\Gamma$, where $G$ is a nilpotent Lie groups of step at most $2$, and $\Gamma$ is a discrete cocompact subgroup.  A \emph{fundamental $2$-step nilmanifold} is one of the following three examples of a $2$-step nilmanifold:
\begin{itemize}
\item \textup{(Unit circle)} $G = \R$ and $\Gamma = \Z$.
\item \textup{(Skew torus)} $G = \left(\begin{smallmatrix} 1 & \Z & \R \\ 0 & 1 & \R \\ 0 & 0 & 1 \end{smallmatrix}\right)$ and $\Gamma = \left(\begin{smallmatrix} 1 & \Z & \Z \\ 0 & 1 & \Z \\ 0 & 0 & 1 \end{smallmatrix}\right)$.
\item \textup{(Heisenberg nilmanifold)} $G = \left(\begin{smallmatrix} 1 & \R & \R \\ 0 & 1 & \R \\ 0 & 0 & 1 \end{smallmatrix}\right)$ and $\Gamma = \left(\begin{smallmatrix} 1 & \Z & \Z \\ 0 & 1 & \Z \\ 0 & 0 & 1 \end{smallmatrix}\right)$.
\end{itemize}
We place smooth metrics on each of these nilmanifolds; the exact choice of metric is not important. An \emph{elementary $2$-step nilmanifold} is a Cartesian product of finitely many fundamental $2$-step nilmanifolds, with the product metric. Again, the exact convention for defining the product metric is not important.  An \emph{elementary $2$-step nilsequence} is a sequence of the form $n \mapsto F(g^n x_0)$, where $G/\Gamma$ is an elementary $2$-step nilmanifold, $F: G/\Gamma \to \C$ is a Lipschitz function, $g \in G$, and $x_0 \in \Gamma$.
\end{definition}

\emph{Remarks.}  If one only had the unit circle and not the skew torus and Heisenberg nilmanifold, the notion of an elementary $2$-step nilsequence would collapse to that of a quasiperiodic sequence.  It is not hard to see that the unit circle and skew torus can be embedded into the Heisenberg nilmanifold, and so one may work entirely with products of Heisenberg nilmanifolds if one wished.  For further discussion of nilsequences see \cite{bhk,bl,green-tao-u3inverse,green-tao-nilratner,hkanal}.

\begin{theorem}[Inverse theorem for $U^3(\Z/N\Z)$]\label{invu3}
Suppose that $f : \Z/N\Z \rightarrow \C$ is a $K$-approximate quadratic for some $K \geq 2$. Then $f$ $\exp(-K^C)$-correlates with an elementary $2$-step nilsequence $F(g^n x_0)$, where $F: G/\Gamma \to \C$ is Lipschitz with Lipschitz constant at most $\exp(K^C)$, $g \in G$, $x_0 \in G/\Gamma$ and $G/\Gamma$ is an elementary $2$-step nilsystem of dimension at most $K^C$.
\end{theorem}

\emph{Remarks.} The proofs of Theorems \ref{invu3-f2} and \ref{invu3} depend very heavily on earlier work of Gowers \cite{gowers:ap4,gowers-longaps}. In Theorem \ref{invu3} one can replace the notion of an elementary $2$-step nilsequence $n \mapsto F(g^n x_0)$ by the more concrete notion of a \emph{bracket phase polynomial}
\begin{equation}\label{brack}
n \mapsto e( \sum_{j=1}^d \alpha_j \{ \beta_j n \} \{ \gamma_j n \} + \sum_{k=1}^{d'} \delta_k \{ \eta_k n \} )
\end{equation}
where $\alpha_j, \beta_j, \gamma_j, \delta_k, \eta_k \in \R$, $\{x\}$ is the fractional part of $x$ (defined to lie in $(-1/2,1/2]$), and $d, d'$ are integers of size at most $K^C$. See \cite{bl,green-tao-u3inverse} for further discussion.\vspace{11pt}

Once again, it is not generally thought that the bounds in these two results are best possible. The following two conjectures might be referred to as the \emph{Polynomial inverse conjectures for the $U^3$ Gowers norms}, or $\PGI(3)$ for short.

\begin{conjecture}[$\PGI(3)$ over $\F_2^{n}$]\label{PGI-1}
Suppose that $f : \F_2^n \rightarrow \C$ is a $K$-approximate quadratic. Then $f$ $K^{-C}$-correlates with a quadratic phase $(-1)^{\psi}$.
\end{conjecture} 

\begin{conjecture}[Weak $\PGI(3)$ over $\Z/N\Z$]\label{PGI-2}
Suppose that $f : \Z/N\Z \rightarrow \C$ is a $K$-approximate quadratic. Then $f$ $\exp(-K^{o(1)})$-correlates with an elementary $2$-step nilsequence $F(g^n x_0)$, where $F: G/\Gamma \to \C$ is Lipschitz of order at most $\exp(K^{o(1)})$, $g \in G$, $x_0 \in G/\Gamma$ and $G/\Gamma$ is an elementary $2$-step nilsystem of of dimension at most $K^{o(1)}$.
\end{conjecture}

\emph{Remarks.} The second of these conjectures deserves some comment. Usually, when inverse conjectures for the Gowers norms are discussed (for example in \cite{green-tao-linearprimes}) there is no restriction to \emph{elementary} nilsequences. We have made this restriction here to simplify the discussion, and in particular to avoid the need to involve the quantitative theory of 2-step nilmanifolds in general as was done in the first two sections of \cite{green-tao-nilratner}. However it transpires that Conjecture \ref{PGI-2} is implied by the same conjecture without the restriction to elementary nilsequences, simply because every $2$-step nilsequence may be closely approximated by a weighted sum of elementary 2-step nilsequences. We omit the details of this deduction, which can be obtained from the calculations in Appendix B of \cite{green-tao-u3mobius}.

We do not dare, at this stage, to even formulate a strong $\PGI(3)$ conjecture over $\Z/N\Z$. To do so would appear to involve rather subtle issues connected with the exact definition of complexity of a nilsequence.\vspace{11pt}
 
We are now in a position to state our main results.
 
\begin{theorem}[Equivalence of PFR and PGI(3), finite field version]\label{mainthm-1}
 Conjecture \ref{PFR-2} and Conjecture \ref{PGI-1} are equivalent.
\end{theorem}
 
\emph{Remark.} A similar result would hold over $\F_p$, for any fixed prime $p$, though the exponents obtained would depend on $p$. 
 
 \begin{theorem}[Equivalence of PFR and PGI(3), $\Z$-version]\label{mainthm-2}
 Conjecture \ref{PFR-1} and Conjecture \ref{PGI-2} are equivalent.
 \end{theorem}

 The fact that Conjecture \ref{PFR-2} implies Conjecture \ref{PGI-1} follows by a modification of Samorodnitsky's argument \cite{samorodnitsky}, and similarly the fact that Conjecture \ref{PFR-1} implies Conjecture \ref{PGI-2} follows from modification of \cite{green-tao-u3inverse}.  Both arguments are strongly dependent on the work of Gowers mentioned earlier. The details of these deductions are a little technical and are discussed in Appendix \ref{gowi}.  However, the main novelty of our paper lies in the opposite implications $\PGI(3) \Rightarrow \PFR$, the discussion of which forms the main body of this paper.\vspace{11pt}
 
\emph{Remark.} The methods used to prove Theorems \ref{mainthm-1}, \ref{mainthm-2} also establish an equivalence between Theorem \ref{inv1} and Theorem \ref{invu3-f2}, and between Theorem \ref{inv2} and Theorem \ref{invu3}, though such an equivalence is redundant given that all four theorems have already been proven in the literature.\vspace{11pt}
 
Let us conclude by remarking that Shachar Lovett informed us that he independently observed Theorem \ref{mainthm-1}.
 
\section{The finite field case}\label{ffc}

We turn now to the proof that Conjecture \ref{PFR-2} implies Conjecture \ref{PGI-1}, that is to say the $\PGI(3)$ implies the $\PFR$ over the finite field $\F_2$.  The argument proceeds via the following intermediate result concerning the structure of \emph{approximate homomorphisms} on the infinite vector space $\F_2^\infty := \bigcup_n \F_2^n$.

\begin{lemma}[Approximate homomorphisms]\label{approx-hom-f2}
Assume Conjecture \ref{PGI-1}. Suppose that $S \subseteq \F_2^n$ is a set of cardinality $\sigma 2^n$ for some $0 < \sigma < 1/2$, and that $\phi : S \rightarrow \F_2^{\infty}$ is a \emph{Fre\u{\i}man homomorphism} on $S$, i.e. $\phi(x_1)+\phi(x_2) = \phi(x_3)+\phi(x_4)$ whenever $x_1,x_2,x_3,x_4 \in S$ are such that $x_1+x_2=x_3+x_4$. Then there is an affine linear map $\psi : \F_2^n \rightarrow \F_2^{\infty}$ such that $\phi(x) = \psi(x)$ for at least $\sigma^C 2^n$ values of $x \in S$.
\end{lemma}

\emph{Remark.} By combining this lemma with known additive-combinatorial results one could obtain the conclusion of Lemma \ref{approx-hom-f2} under \emph{a priori} weaker assumptions, for example that $\P(\phi(x_1) + \phi(x_2) = \phi(x_3) + \phi(x_4) | x_1 + x_2 = x_3 + x_4)$ is large. Indeed a map of this type restricts to a Fre\u{\i}man homomorphism on a large set $S$ by arguments of Gowers and Ruzsa (see \cite[Section 7]{gowers-longaps}).\vspace{11pt}

Let us first show how Conjecture \ref{PFR-1} follows from Lemma \ref{approx-hom-f2}. Suppose that $A \subseteq \F_2^{\infty}$ is a $K$-approximate group, and let $n$ be minimal such that there exists a linear map $\pi : \F_2^{\infty} \rightarrow \F_2^n$ which is a Fre\u{\i}man isomorphism\footnote{A \emph{Fre\u{\i}man isomorphism} is a Fre\u{\i}man homomorphism which is invertible and whose inverse is also a Fre\u{\i}man homomorphism.} from $A$ to $\pi(A)$; this quantity $n$, which one can view as a sort of ``rank'' or ``dimension'' for $A$, is finite since $A$ is finite. If there is some\footnote{We use $kA = A + \ldots + A$ to denote the $k$-fold iterated sumset of $A$, thus $4\pi(A) = \pi(A)+\pi(A)+\pi(A)+\pi(A)$.  Note in $\F_2$ that there is no distinction between sums and differences, thus for instance $4\pi(A) = 2\pi(A)-2\pi(A)$.} $x \in\F_2^n \setminus 4\pi(A)$ then we could compose $\pi$ with the projection map $\psi : \F_2^n \rightarrow \F_2^n /\langle x\rangle$ to obtain a linear map $\pi : \F_2^{\infty} \rightarrow \F_2^{n-1}$ which is a Fre\u{\i}man isomorphism when restricted to $A$, contrary to the assumed minimality of $n$. It follows that $4\pi(A) = \F_2^n$. But $\pi(A)$ is Fre\u{\i}man isomorphic to $A$, which is a $K$-approximate group. It follows that the doubling constant $|2 \pi(A)|/|\pi(A)|$ is at most $K$, and hence by Ruzsa's sumset estimates (cf. \cite[Corollary 2.23]{tao-vu}) that $2^n = |4\pi(A)| \leq K^C|A|$.

What we have done here is find a ``dense model'' $\pi(A) \subset \F_2^n$ of the set $A$; the simple argument we used to do so is the finite field analogue of an argument of Ruzsa \cite{ruzsa-freiman} that we shall recall later in the paper. Write $S = \pi(A)$ and $\phi$ for the inverse of $\pi$, restricted to $S$. Then $\phi$ is a Fre\u{\i}man homomorphism on $S$ and the set $A$ is precisely the image $\phi(S)$. Applying Lemma \ref{approx-hom-f2}, we see that at least $K^{-C}|A|$ of the elements of $A$ are contained in a coset of an $n$-dimensional subspace $H \leq \F_2^{\infty}$. Finally, it follows immediately from standard covering lemmas (cf. \cite[Section 2.4]{tao-vu}) that $A$ is $K^C$-controlled by $H$. \endproof

It remains, then, to establish Lemma \ref{approx-hom-f2}.  The key observation linking Fre\u{\i}man homomorphisms to approximate quadratics is the following lemma.

\begin{lemma}\label{construct-lem} Suppose that $S \subseteq \F_2^n$ is a set of size $\sigma 2^n$ for some $0 < \sigma < 1/2$ and that $\phi : S \rightarrow \F_2^{\infty}$ is a Fre\u{\i}man homomorphism. The image of $\phi$ certainly lies in some finite-dimensional subspace $\F_2^N$.  If $f : \F_2^{n+N} \rightarrow [-1,1]$ is the function $f(x,y) := 1_S(x)(-1)^{\phi(x) \cdot y}$, then
$\Vert f \Vert_{U^3(\F_2^{n+N})} \geq \sigma$.
\end{lemma}

\proof Consider a parallelopiped $(x+\omega \cdot h, y + \omega \cdot k)_{\omega \in \{0,1\}^3}$ in the support of $f$, where $h=(h_1,h_2,h_3)$, $k = (k_1,k_2,k_3)$, $x,h_1,h_2,h_3 \in \F_2^n$ and $y,k_1,k_2,k_3 \in \F_2^N$.  Then $x + \omega \cdot h \in S$ for all $\omega \in \{0,1\}^3$.  Since $\phi$ is a Fre\u{\i}man homomorphism on $S$, we see that $\phi(x+\omega \cdot h)$ depends linearly on $\omega$, and so $\phi(x+\omega \cdot h) \cdot (y + \omega \cdot k)$ depends quadratically on $\omega$.  Since $\{0,1\}^3$ is three-dimensional, we conclude that
$$ \sum_{\omega \in \{0,1\}^3} (-1)^{|\omega|} \phi(x+\omega \cdot h) \cdot (y + \omega \cdot k) = 0$$
where $|\omega|$ is the number of $1$s in the coefficients of $\omega$ (actually, as we are working in $\F_2$ here, the $(-1)^{|\omega|}$ factor could in fact be ignored).  From this and the definition of $f$ and the $U^3(\F_2^{n+N})$ norm we conclude that
$$ \| f \|_{U^3(\F_2^{n+N})} = \| 1_S \|_{U^3(\F_2^{n+N})}.$$
The behaviour in the $y$ index is now trivial, and therefore
$$ \| 1_S \|_{U^3(\F_2^{n+N})} = \| 1_S \|_{U^3(\F_2^{n})}.$$
Meanwhile, $\|1_S\|_{U^1(\F_2^n)} \geq \sigma$ by hypothesis.  The claim now follows from the monotonicity of the Gowers norms (see, for example, \cite[eq. 11.7]{tao-vu}).
\endproof\vspace{11pt}

Now suppose $S, \sigma$ are as in the statement of Lemma \ref{approx-hom-f2}, and let $N$ and $f$ be as in the above lemma. 
Assuming Conjecture \ref{PGI-1} for this choice of $f$, there exists a quadratic polynomial $\Psi: \F_2^{n+N} \to \F_2$ such that
$$ |\E_{x \in \F_2^n} \E_{y \in \F_2^N} 1_S(x) (-1)^{\phi(x) \cdot y} (-1)^{\Psi(x,y)}| \geq \sigma^C.$$
Thus, for at least $\geq \sigma^C 2^n$ values of $x \in S$, one has
\begin{equation}\label{yo}
 |\E_{y \in \F_2^N} (-1)^{\phi(x) \cdot y} (-1)^{\Psi(x,y)}| \geq \sigma^C.
\end{equation}
Let us fix $x$ so that \eqref{yo} holds.  We may split $\Psi(x,y)$ as 
\begin{equation}\label{psid}
 \Psi(x,y) = \Psi(0,y) + \Psi(x,0) - \Psi(0,0) + B(x,y)
 \end{equation}
where $B$ is the ``mixed derivative'' of $\Psi$, defined as
$$ B(x,y) := \Psi(x,y) - \Psi(x,0) - \Psi(0,y) + \Psi(0,0).$$
From \eqref{yo} it thus follows that 
$$ |\E_{y \in \F_2^N} (-1)^{\phi(x) \cdot y} (-1)^{B(x,y)} (-1)^{\Psi(0,y)}| \geq \sigma^C.$$
As $\Psi$ is quadratic, $B$ is bilinear in $x$ and $y$, and hence $B(x,y) = \psi(x) \cdot y$ for some linear map $\psi: \F_2^N \to \F_2^n$.  We conclude that
$$ |\E_{y \in \F_2^N} (-1)^{(\phi(x)-\psi(x)) \cdot y} (-1)^{\Psi(0,y)}| \geq \sigma^C,$$
which means that the function $y \mapsto (-1)^{\Psi(0,y)}$ has a Fourier coefficient of size at least $\sigma^C$ at $\phi(x)-\psi(x)$.  Hence by Plancherel's theorem the number of such large Fourier coefficients is at most $\sigma^{-2C}$.  We conclude that $\phi(x)-\psi(x)$ takes at most $\sigma^{-2C}$ values on at least $\sigma^C 2^n$ values of $x \in S$, and the claim follows from the pigeonhole principle.\endproof 

\section{The integer case}

We turn now to the proof that Conjecture \ref{PGI-2} implies Conjecture \ref{PFR-2}. This argument goes along similar lines to that in the previous section, but is somewhat more involved since one must deal with nilsequences rather than quadratic forms. We present the argument in such a way as to emphasise the close parallels with the preceding section.

Once again matters rest on a reduction to an inverse theorem for approximate homomorphisms. We write $[N]$ for the set $\{1,\ldots,N\}$.

\begin{lemma}[Approximate homomorphisms]\label{approx-hom-z}
Assume Conjecture \ref{PGI-2}. Suppose that $N$ is a positive integer, that $S \subseteq [N]$ is a set of cardinality $\sigma N$, and that $\phi : S \rightarrow \Z$ is a Fre\u{\i}man homomorphism on $S$. Then there is a generalised arithmetic progression $P \subseteq [N]$ of dimension $\sigma^{-o(1)}$ and size at least $\exp(-\sigma^{-o(1)})N$ together with a Fre\u{\i}man homomorphism $\psi : P \rightarrow \Z$ such that $\phi(x) = \psi(x)$ for at least $\exp(-\sigma^{-o(1)})N$ values of $x \in S$.
\end{lemma}

The proof that this lemma implies Conjecture \ref{PFR-2} is not particularly onerous and goes along much the same lines as the argument in the previous section. Supposing that $A \subseteq \Z$ is a $K$-approximate subgroup, Ruzsa's ``model lemma'' \cite[Theorem 2]{ruzsa-model} implies that there is a $N \leq K^C|A|$ together with a subset $A' \subseteq A$ of cardinality at least $|A|/2$ and a Fre\u{\i}man isomorphism $\pi : A' \rightarrow S$ to a subset $S \subseteq [N]$. Write $\phi  := \pi^{-1}$, and observe that $\phi: S \to \Z$ has image $\phi(S)=A'$.  Noting that $|S| \geq K^{-C}N$, it follows from Lemma \ref{approx-hom-z} and the fact that Fre\u{\i}man isomorphisms take generalised progressions to generalised progressions (see \cite[Proposition 5.24]{tao-vu}) that at least $\exp(-K^{o(1)})|A|$ of $A$ is contained in a generalised progression in $\Z$ of dimension $K^{o(1)}$ and cardinality at most $N \leq K^C|A|$. Once again, standard covering arguments complete the deduction of Conjecture \ref{PFR-2}.\endproof

It remains to prove Lemma \ref{approx-hom-z}.  The starting point is the following analogue of Lemma \ref{construct-lem}, showing how to convert Fre\u{\i}man homomorphisms to approximate quadratics.

\begin{lemma}\label{construct-z} Let $N, M \geq 1$ be integers, let $S \subset [N]$ be such that $|S| \geq \sigma N$, and let $\phi: S \to \Z/M\Z$ be a Fre\u{\i}man homomorphism.  Define a function $f : \Z/4NM\Z \rightarrow \C$ by
\[ f(x + 4Ny) := \left\{ \begin{array}{ll} 1_S(x)e_M(\phi(x)y) & \mbox{if $x \in [N], y \in \Z/M\Z$}; \\ 0 & \mbox{otherwise},  \end{array}  \right.\]  where $e_M(x) := e^{2\pi ix/M}$, and $4Ny \in \Z/4NM\Z$ is defined in the obvious manner for $y \in \Z/M\Z$. Then $\Vert f \Vert_{U^3(\Z/4 NM\Z)} \geq \frac{1}{4} \sigma$.
\end{lemma}

\proof Every parallelepiped in the support of $f$ is of the form $(x+\omega \cdot h, 4N(y + \omega \cdot k))_{\omega \in \{0,1\}^3}$, where $y \in \Z/M\Z$, $k = (k_1,k_2,k_3) \in (\Z/M\Z)^3$, and $x + \omega \cdot h \in S$.  By arguing exactly as in Lemma \ref{construct-lem} we have that
$$ \sum_{\omega \in \{0,1\}^3} (-1)^{|\omega|} \phi( x + \omega \cdot h ) ( y + \omega \cdot k ) = 0$$
and so
$$ \|f\|_{U^3(\Z/4NM\Z)} = \| 1_{\tilde S} \|_{U^3(\Z/4NM\Z)}$$
where $\tilde S := \{ x+4Ny: x \in S; y \in \Z/M\Z\}$ is the support of $f$.  But we have
$$ \| 1_{\tilde S} \|_{U^1(\Z/4NM\Z)} \geq \sigma/4$$
and the claim follows as before from the monotonicity of the Gowers norms.
\endproof\vspace{11pt}

We return now to the proof of Lemma \ref{approx-hom-z}.  
That lemma deals with Fre\u{\i}man homomorphisms $\phi :S \rightarrow \Z$. However such a map is a Fre\u{\i}man homomorphism  if and only if the composition $\pi_M \circ \phi$ is a Fre\u{\i}man homomorphism for all sufficiently large $M$, and so we may suppose instead that $\phi$ maps $S$ to $\Z/M\Z$ for some $M$. 

Let $f$ be as in Lemma \ref{construct-z}.
Assuming Conjecture \ref{PGI-2}, it follows that there is some elementary $2$-step nilmanifold $G/\Gamma$ of dimension at most $\sigma^{-o(1)}$, a function $F: G/\Gamma \to \C$ of Lipschitz constant at most $\exp(\sigma^{-o(1)})$, $g \in G$, and $x_0 \in G/\Gamma$ such that
$$
|\E_{x \in [N]}\E_{y \in [M]} 1_S(x) e_M(\phi(x)y) F(g^{x + 4Ny} x_0)| \geq \exp(-\sigma^{-o(1)}).
$$
Writing $x_0 = g_0 \Gamma$ for some $g_0 \in G$ of distance at most $\exp(\sigma^{-o(1)})$ from the origin, and rewriting $g^{x+4Ny}x_0 = g_0 \tilde g^{x+4Ny} \Gamma$ where $\tilde g := g_0^{-1} g g_0$, we see (after shifting $F$ by $g_0$ and replacing $g$ by $\tilde g$ if necessary) that we may normalise $x_0$ to be at the origin $\Gamma$.  By embedding the skew torus in the Heisenberg group if necessary we may take $G$ to be a product of Heisenberg groups and hence, in particular, connected and simply-connected.

The vertical torus $[G,G] / (\Gamma \cap [G,G])$ of the elementary $2$-step nilmanifold can be identified with a torus $(\R/\Z)^{d_2}$ for some $d_2 \leq \sigma^{-o(1)}$. By standard harmonic analysis arguments (see, for example, \cite[Lemma A.9]{green-tao-u3mobius}) the Lipschitz function may be decomposed into a linear combination of at most $\exp(\sigma^{-o(1)})$ Fourier characters along the vertical direction with coefficients of magnitude at most $\exp(\sigma^{-o(1)})$, plus an error of $\exp(-\sigma^{-o(1)})$ in $L^{\infty}$. Applying the pigeonhole principle it follows that one may assume that $F$ is a \emph{vertical character}, which means that there exists a character $\chi: [G,G] / (\Gamma \cap [G,G]) \to S^1$ such that
\begin{equation}\label{fosc}
 F( g_2 x ) = \chi(g_2) F(x)
\end{equation}
for all $x \in G/\Gamma$ and $g_2 \in [G,G]$ (where we lift $\chi$ to $[G,G]$ in the obvious fashion).  

The Lipschitz function $|F|$ is now invariant under the action of the vertical torus and descends to a function on the \emph{horizontal torus} $G /\Gamma [G,G]$, which can be identified with a torus $(\R/\Z)^{d_1}$ for some $d_1 \leq \sigma^{-o(1)})$.  By applying a Lipschitz partition of unity we may assume that $|F|$ (and hence $F$) is supported in a small ball in this torus, of radius less than $\exp(-\sigma^{-o(1)})$ say.

By the pigeonhole principle, we can now find $\geq \exp(-\sigma^{-o(1)})N$ values of $x \in S$ such that
$$
|\E_{y \in [M]} e_M(\phi(x)y) F(g^{x + 4Ny} \Gamma)| \geq \exp(-\sigma^{-o(1)}).
$$

By pigeonholing in $x$ (reducing the number of available $x$ by a factor of $\exp(-\sigma^{-o(1)})$), we may assume that for all these $x$ the point $g^x \Gamma$ lies in a small ball $B$ in $G/\Gamma$, of radius less than $\exp(\sigma^{-o(1)})$.

We turn now to the task of simplifying $F(g^{x + 4Ny} \Gamma)$: this may be thought of, roughly, as a quest to find a suitable analogue for the decomposition \eqref{psid}. To begin with let us expand $g^x$ as $\{g^x\} \lfloor g^x\rfloor$, where $\lfloor g^x\rfloor \in \Gamma$ and $\{g^x\}$ lies in a fundamental domain of $G/\Gamma$ that contains $B$ in its interior\footnote{Several papers of the authors -- for example the appendix of \cite{green-tao-u3mobius} -- contain example computations of $\{g^x\}$ and $\lfloor g^x\rfloor$ on the Heisenberg group for fundamental domains like $\{-\frac{1}{2},\frac{1}{2}\}$ or $[0,1]^3$.}. As usual, write $[g,h] := ghg^{-1}h^{-1}$ for the commutator of two elements $g$ and $h$ in some ambient group. Now in any 2-step nilpotent group $G$ we have $[x^n,y] = [x,y]^n$ for all $x,y \in G$ and all $n \in \Z$: this follows from the commutator identity $[xy,z] = [y,z]^x[x,z]$, which is valid in all groups. It follows that
$$ g^{x+4Ny} \Gamma = g^{4Ny} \{g^x\} \Gamma = [g^{4N},\{g^x\}]^y \{g^x\} g^{4Ny} \Gamma.$$
 Since $F$ is a vertical character, we thus see that
$$ F(g^{x+4Ny} \Gamma) = \chi([g^{4N},\{g^x\}])^y F( \{g^x\} g^{4Ny} \Gamma )$$
and so
$$
|\E_{y \in [M]} e([\frac{1}{M} \phi(x)-\psi(x)]y) F(\{g^x\} g^{4Ny} \Gamma)| \geq \exp(-\sigma^{-o(1)})
$$
for at least $\exp(-\sigma^{-o(1)})N$ values of $x \in S$, where $\psi(x) \in \R/\Z$ is the phase such that
$$\chi([g^{4N},\{g^x\}]) = e(\psi(x)).$$

By construction, $\{g^x\}$ is supported in a small ball centred at some $g_0 \in G$, of radius less than $\exp(-\sigma^{-o(1)})$. Provided that this ball is chosen small enough, the Lipschitz nature of $F$ guarantees that
$$
|\E_{y \in [M]} e([\frac{1}{M} \phi(x)-\psi(x)]y) F(g_0 g^{4Ny} \Gamma)| \geq \exp(-\sigma^{-o(1)}).
$$
Recall that $|F|$ has small support, on account of the partition of unity that was brought into play earlier in the argument. We now let $F_0: G/\Gamma \to \C$ be another function of Lipschitz constant at most $\exp(\sigma^{-o(1)})$ and with vertical character $\chi$ which has magnitude $1$ on the support of $F(g_0 \cdot)$; there are no topological obstructions to building such an $F_0$ if the support of $|F|$ is small enough (think, for example, of the function $\psi(x,y,z)e(z)$ on the Heisenberg nilmanifold $G/\Gamma$, where $\psi$ is supported on a small ball and equals $1$ on a very small ball about the origin in the fundamental domain $\{-\frac{1}{2},\frac{1}{2}\}$).

With this function $F_0$ constructed we may write
$$ F(g_0 g^{4Ny} \Gamma) = \tilde F( g^{4Ny} \Gamma ) F_0( g^{4Ny} \Gamma)$$
where $\tilde F: G/\Gamma \to \C$ is the function
$$ \tilde F( x ) := F( g_0 x ) \overline{F_0}(x).$$
Observe that the function $\tilde F(x)$ is invariant under the action of the vertical torus, and thus descends to a function on $(\R/\Z)^{d_1}$, which by abuse of notation we also call $\tilde F$. Thus
$$ F(g_0 g^{4Ny} \Gamma) = \tilde F( \pi(g^{4Ny} \Gamma) ) F_0( g^{4Ny} \Gamma),$$
where $\pi: G/\Gamma \to (\R/\Z)^{d_1}$ is the projection onto the horizontal torus.

The function $\tilde F$ is Lipschitz with constant at most $\exp(\sigma^{-o(1)})$, and so (by \cite[Lemma A.9]{green-tao-u3mobius}) can be decomposed into a combination of at most $\exp(\sigma^{-o(1)})$ characters with coefficients at most $\exp(\sigma^{-o(1)})$, plus an error of size $\exp(-\sigma^{-o(1)})$.  Meanwhile, $\pi(g^{4Ny} \Gamma) \in (\R/\Z)^{d_1}$ evolves linearly in $y$.  By the pigeonhole principle, refining the set of available $x$ some more, we may thus assume that
\begin{equation}\label{to-apply-ls}
|\E_{y \in [M]} e([\frac{1}{M} \phi(x)-\psi(x)]y) e( \xi_0 y ) F_0(g^{4Ny} \Gamma)| \geq \exp(-\sigma^{-o(1)}).
\end{equation}
for some $\xi_0 \in \R/\Z$ independent of $x$, and for at least $\exp(-\sigma^{-o(1)}) N$ values of $x$.  Thus, the function $y \mapsto F_0(g^{4Ny} \Gamma)$ has a large Fourier coefficient at $\frac{1}{M} \phi(x)-\psi(x)+\xi_0$.  

In the finite field argument we applied Plancherel's theorem at this point. Here the appropriate tool is the large sieve, a kind of approximate version of Plancherel which states that a function $f : [M] \rightarrow \C$ cannot have large Fourier coefficients at many \emph{separated} points. The following (standard) statement of it may be found in \cite[Ch. 27]{davenport}: if the points $\theta_1,\dots,\theta_K \in \R/\Z$ are $\delta$-separated then 
\[ \sum_{j=1}^K |\sum_{y \in [M]} f(y) e(y\theta_j)|^2 \ll (M + \delta^{-1})\sum_{y \in [M]} |f(y)|^2.\]
Applying this to \eqref{to-apply-ls} and the remark following it, we see that the large Fourier coefficients $\frac{1}{M} \phi(x)-\psi(x)+\xi_0$ of the function $y \mapsto F_0(g^{4Ny}\Gamma)$ can be covered by at most $\exp(\sigma^{-o(1)})$ arcs of length $1/M$ on the unit circle $\R/\Z$. Pigeonholing, and refining the set of $x$ by yet another factor of $\exp(-\sigma^{-o(1)})$, we may assume that $\frac{1}{M} \phi(x)-\psi(x)+\xi_0$ lies inside a fixed arc of length $\frac{1}{100M}$.  This implies, refining the set of $x$ one more time, that we may find a $\xi_1 \in \R/\Z$ such that for at least $\exp(-\sigma^{-o(1)}) N$ values of $x \in S$, $\frac{1}{M} \phi(x)-\psi(x)+\xi_1 \in [-1/100M,1/100M]$.

By direct computations on the Heisenberg group along the lines of those in \cite{green-tao-u3mobius} we see that $\pi(\{g^x\}) = (\alpha_1 x,\ldots,\alpha_{d_1} x)$ for some $\alpha_1,\ldots,\alpha_{d_1} \in \R/\Z$, and then that $$ \chi([g^{4N},\{g^x\}]) = e( \sum_{j=1}^{d_1} \beta_j \{ \alpha_j x - \gamma_j \} )$$
for some $\beta_j, \gamma_j \in \R/\Z$ independent of $x$. Here the fractional part $\{t\}$ of $t \in \R$ is chosen to lie in $(-\frac{1}{2},\frac{1}{2}]$,  and the need for the shift $\gamma_j$ arises from the fact that $\{g^x\}$ is chosen to lie in a fundamental domain of $G/\Gamma$ containing $B$ in its interior.

This means, of course, that
$$ \psi(x) = \sum_{j=1}^{d_1} \beta_j \{ \alpha_j x - \gamma_j\}.$$

The set of all $x \in [N]$ such that $\pi(g^x \Gamma)$ lies within $\exp(-\sigma^{-o(1)})$ of the origin is a Bohr set of rank at most $\sigma^{-o(1)}$ and radius at least $\exp(-\sigma^{-o(1)})$, and hence by \cite[Theorem 3.1]{ruzsa-freiman} (reproduced as Theorem \ref{bohr-prog} in the appendix) it contains a proper symmetric generalised arithmetic progression $P$ of dimension at most $\sigma^{-o(1)}$ and cardinality at least $\exp(-\sigma^{-o(1)}) N$. By discarding generators of $P$ if necessary we may assume that all sidelengths of $P$ are at least $C_0$ for some constant $C_0$ to be specified later.  By standard covering lemmas such as \cite[Lemma 2.14]{tao-vu} we may cover $[N]$ by at most $\exp(\sigma^{-o(1)} N)$ translates of $P$, so by the pigeonhole principle we may assume that all the $x$ under discussion, that is to say those $x$ for which $\alpha_j x \approx \gamma_j$, are contained in a single translate $x_0+P$ of $P$. Note that each map $x \mapsto \{\alpha_j x - \gamma_j\}$ is a Freiman homomorphism on $x_0 + P$ and hence so is the entire phase $\psi$.

If we let $Q$ be the set of all $x \in x_0+P$ such that $\psi(x)-\xi_1$ lies within $\frac{1}{100 M}$ of a multiple $\frac{1}{M} \tilde \phi(x)$ of $\frac{1}{M}$, where $\tilde \phi(x) \in \Z/M\Z$, then we conclude upon rounding to the nearest multiple of $\frac{1}{M}$ that $\tilde \phi$ is a Freiman homomorphism on $Q$.  Also, from construction we see that $\phi(x) = \tilde \phi(x)$ for at least $\exp( - \sigma^{-o(1)}) N$ values of $x \in S \cap Q$.

To conclude the argument one needs to show that $Q$ contains a generalised arithmetic progression of dimension at most $\sigma^{-o(1)}$ and cardinality at least $\exp(-\sigma^{-o(1)}) N$ (since one can then cover $Q$ by at most $\exp(\sigma^{-o(1)})$ translates of such a progression). This will follow straightforwardly from the following lemma which, though it looks to be of a standard type, does not appear to be in the literature. A proof may be found in Appendix \ref{prog-app}.

\begin{a1} Let $\eps \in (0,1/2)$ be a real number.
Suppose that $P$ is a $d$-dimensional proper progression with sidelengths $N_1,\dots,N_d > C/\eps$ and that $\eta : P \rightarrow \R/\Z$ is a Fre\u{\i}man homomorphism which vanishes at some point of $P$. Then the set $\{x \in P : \Vert \eta(x)\Vert_{\R/\Z} \leq \eps\}$ contains a progression of dimension at most $d+1$ and size at least $(Cd)^{-d}\eps^{d+1}|P|$.
\end{a1}

We shall apply the lemma with $\eps = 1/100$, this being valid if the constant $C_0$ was chosen to be large enough earlier on. Recall that there are many $x \in S \cap (x_0 + P)$ such that $\frac{1}{M}\phi(x) - \psi(x) + \xi_1 \in [-1/100M,1/100M]$. Pick one such $x^*$, and take $\xi_2$ to be such that $\frac{1}{M}\phi(x^*) - \psi(x^*) + \xi_2 = 0$ and $\Vert \xi_1 - \xi_2\Vert_{\R/\Z} \leq 1/100M$.  Now we simply apply Lemma \ref{prog-lem} to the progression $x_0 + P$, taking $\eta = M(\psi - \xi_2)$ and $\eps = 1/100$.  The progression $Q$ contains the set $\{x \in x_0 + P : \Vert \eta(x)\Vert_{\R/\Z} \leq 1/100\}$, and of course $\eta$ vanishes at $x^*$. It follows from Lemma \ref{prog-lem} that $Q$ does indeed contain a generalised arithmetic progression of dimension at most $\sigma^{-o(1)}$ and cardinality at least $\exp(-\sigma^{-o(1)}) N$, and this concludes the proof of Theorem \ref{mainthm-2}.\endproof

\section{Higher order correspondences}

It appears that the correspondence between inverse sumset theorems and inverse conjectures for the Gowers norms have some partial higher order analogues, although the situation here is much less well understood.  To illustrate this phenomenon, consider the following result, recently proven in \cite{btz,tz}. Here and for the rest of the section we write $\F := \F_5$ for definiteness, although the same arguments would work for $\F_p$ for any fixed prime $p \geq 5$. There are definite issues in extremely low characteristic: see for example \cite{gt-ff-ratner,lms}.

\begin{theorem}[GI(4) over $\F^n$]\label{conj-4} For every $K \geq 2$ there exists an $\eps > 0$ such that if $f : \F^n \rightarrow \C$ is a $K$-approximate cubic in the sense that $\|f\|_\infty \leq 1$ and $\|f\|_{U^4(\F^n)} \geq 1/K$, then $f$ $\eps$-correlates with a \textup{(}genuine\textup{)} cubic phase $e_\F(\psi)$, where $e_{\F}(x) := e^{2\pi i x/|\F|}$ and $\psi: \F^n \to \F$ is cubic in the sense that $\Delta_{h_1} \ldots \Delta_{h_4} \psi(x) = 0$ for all $x,h_1,\ldots,h_4 \in \F^n$.
\end{theorem}

We shall use this theorem to establish the following variant of Lemma \ref{approx-hom-f2}.

\begin{proposition}[Approximate quadratic homomorphisms]\label{approx-hom-f5}
Suppose that $\sigma \in (0,1/2)$, that $S \subseteq \F^n$ is a set of cardinality $\sigma |\F|^n$, and that $\phi : S \rightarrow \F^{\infty}$ is a \emph{Fre\u{\i}man quadratic} on $S$ in the sense that $\sum_{\omega \in \{0,1\}^3} (-1)^{|\omega|} \phi(x + h \cdot \omega) = 0$ whenever $x \in \F^n$ and $h = (h_1,h_2,h_3)$ with $h_1,h_2,h_3 \in \F^n$ are such that $x+\omega \cdot h \in S$.  Then there is a quadratic map $\psi : \F^n \rightarrow \F^{\infty}$ such that $\phi(x) = \psi(x)$ for at least $\eps |\F|^n$ values of $x \in S$, where $\eps = \eps(\sigma) > 0$ depends only on $\sigma$.
\end{proposition}

The initial stages of the proof are very similar to those of Theorem \ref{mainthm-1} and we just sketch them.  As before, we let $N$ be large enough that $\phi$ takes values in $\F^N$, and considers the function $f: \F^{n+N} \to \C$ defined by $f(x,y) := 1_S(x) e_{\F}(\phi(x) \cdot y)$.  A routine modification of Lemma \ref{construct-lem} reveals that
$$ \|f\|_{U^4(\F^{n+N})} \geq \delta$$
and thus by Theorem \ref{conj-4} we can find a cubic $\Psi: \F^{n+N} \to \F$ such that
$$ |\E_{x \in \F^n} \E_{y \in \F^N} 1_S(x) e_{\F}(\phi(x) \cdot y-\Psi(x,y))| \gg_\delta 1,$$
where here we use $X \gg_\delta Y$ to denote the estimate $X \geq C_\delta^{-1} Y$ for some $C_\delta$ depending only on $\delta$.  Thus for $\gg_\delta |\F|^n$ values of $x \in S$, one has
$$ |\E_{y \in \F^N} e_{\F}(\phi(x) \cdot y-\Psi(x,y))| \gg_\delta 1.$$
The next step is to perform a decomposition of $\Psi$ analogous to \eqref{psid}, but unfortunately the analogous decomposition is not so favourable.  Namely, one has
$$ \Psi(x,y) = \Psi(0,y) + Q_x(y) + \psi(x) \cdot y + P(x)$$
where $Q_x: \F^N \to \F$ is a quadratic polynomial that varies affine-linearly in $x$, $\psi: \F^n \to \F^N$ is a quadratic polynomial, and $P: \F^n \to \F$ is a cubic polynomial.  We thus have
\begin{equation}\label{xr}
 |\E_{y \in \F^N} e_{\F}((\phi(x)-\psi(x)) \cdot y -Q_x(y)-\Psi(0,y))| \gg_\delta 1
\end{equation}
for $\gg_\delta |\F|^n$ values of $x \in S$.

The factor of $e_{\F}(-Q_x(y))$ in the functions $f_x(y) := e_{\F}((\phi(x)-\psi(x)) \cdot y-Q_x(y))$ prevents one from immediately using Plancherel's theorem as in Section \ref{ffc}.  However, from standard Gauss sum estimates (see e.g. \cite[Lemma 1.6]{gt-ff-ratner}) we do have
\begin{equation}\label{gauss-est} |\langle f_x, f_{x'} \rangle| \ll |\F|^{-\rk(Q_x-Q_{x'})/2}\end{equation}
for any $x,x'$. Here the \emph{rank} of a quadratic form $Q$ can be defined as the rank of the symmetric matrix describing the homogeneous part of $Q$. By standard linear algebra there is a vector subspace $V_Q \leq \F^n$ with $\dim(V_Q) = \rk(Q)$ such that $Q(y)$ is a quadratic function of the inner products $\langle v,y\rangle$, $v \in V_Q$.

From \eqref{gauss-est} and a standard duality argument related to the large sieve (see, for example, \cite[Ch. 27, Theorem 1]{davenport}) one can show that there cannot exist $k$ different $x_1,\ldots,x_k\in S$ obeying \eqref{xr} with $\rk(Q_{x_i}-Q_{x_j}) \geq k$, if $k$ is large enough depending on $\delta$.  By the greedy algorithm, we may thus find $x_1,\ldots,x_k$ with $k \ll_\delta 1$ such that $\min_{1 \leq i \leq k} \rk(Q_x - Q_{x_i}) \ll_\delta 1$ for all $x$ obeying \eqref{xr}.  By pigeonholing in the $x$ parameter, we conclude that there exists a quadratic form $Q_{x_1}$ such that $\rk(Q_x - Q_{x_1}) \ll_\delta 1$ for $\gg_\delta |\F|^n$ values of $x \in S$.  By translating we may normalise and take $x_1=0$.

Write $Q'_x$ be the homogeneous quadratic component of $Q_x - Q_0$, so that $Q'_x$ depends linearly on $x$ and $\rk(Q'_x) \ll_{\delta} 1$ for $\gg_{\delta} |\F|^n$ values of $x \in S$.
Key to our argument is the following proposition concerning this situation, which may be of independent interest. It states that a linear function to the set of low-rank quadratics must, in a sense, be quite trivial. 

\begin{proposition}[Triviality of linearly varying low-rank quadratic forms] \label{rank-line} Let $r \in \N_0$ and suppose that $\eps \in (0,1]$ is a real number.
Suppose that $x \mapsto Q_x$ is a linear map from $\F^n$ to the space of homogeneous quadratics over $\F^N$. For each such form $Q_x$ associate the vector space $V_x := V_{Q_x}$. Suppose that there is a set $A$ of at least $\alpha |\F|^n$ values of $x$ for which $\rk(Q'_x) \leq r$. Then there is some vector space $V \leq \F^n$, $\dim(V) \leq r$, such that $V_x \subseteq V$ for at least $\alpha'(\alpha,r)|\F|^n$ values of $x \in A$, where $\alpha : (0,1] \times \N_0 \rightarrow \R$ takes positive values.
\end{proposition}
\proof We claim that under the stated hypotheses there is some vector $v$ which lies in at least $\alpha_0(\alpha,r)|\F|^n$ of the spaces $V_x$, where $\alpha_0$ is a function taking positive values. The proposition then follows quickly by induction on $r$, upon passing to a coset of the codimension one subspace $v^{\perp} \leq \F^n$ which contains at least $\alpha |v^{\perp}|$  elements of $A$.

Now by a standard application of Cauchy-Schwarz (see, e.g, \cite[Corollary 2.10]{tao-vu}) there are at least $\alpha^4 |\F|^{3n}$ additive quadruples in $A$, that is to say quadruples $(x_1,x_2,x_3,x_4) \in A^4$ with $x_1 + x_2 = x_3 + x_4$. We say that such a quadruple is \emph{good} if $V_{x_i} \cap (V_{x_j} + V_{x_k}) = \{ 0 \}$ for all 24 choices of distinct $i,j,k \in \{1,2,3,4\}$.

\emph{Case 1.} At least half of the additive quadruples in $A$ are good. Fix a good quadruple $(x_1,x_2,x_3,x_4) \in A^4$. Let $y,h,k \in \F^n$ be arbitrary, and select $h' \in (h + V_{x_1}^{\perp}) \cap (V_{x_2}^{\perp} \cap V_{x_3}^{\perp})$ and $k' \in (k + V_{x_1}^{\perp}) \cap V_{x_4}^{\perp}$. Straightforward linear algebra (and the goodness of the quadruple $(x_1,x_2,x_3,x_4)$) confirms that this is possible.

From the linearity of the map $x \mapsto Q_x$  we have
\[ Q_{x_1} (y) + Q_{x_2} (y) - Q_{x_3}( y) - Q_{x_4} (y) = 0\] and
\[ Q_{x_1} (y + h') + Q_{x_2}(y + h') - Q_{x_3}(y + h') - Q_{x_4} (y + h') = 0.\]
Since $h' \in V_{x_2}^{\perp} \cap V_{x_3}^{\perp}$ the second of these implies that 
\[ Q_{x_1}( y+h') + Q_{x_2}(y) - Q_{x_3}( y) - Q_{x_4}(y + h') = 0.\]
Subtracting the first equation yields
\[Q_{x_1} (y) - Q_{x_1}(y+h') - Q_{x_4} (y) + Q_{x_4}(y+h') = 0.\]
Substituting $y + k'$ for $y$, recalling that $k' \in V_{x_4}^{\perp}$, and subtracting, this implies that
\[ Q_{x_1}(y) - Q_{x_1}(y + h') - Q_{x_1}(y + k') + Q_{x_1}(y + h' + k') = 0.\]
But $h - h'$ and $k - k'$ both lie in $V_{x_1}^{\perp}$, and so this implies that
\[ Q_{x_1}(y) - Q_{x_1}(y + h) - Q_{x_1}(y + k) + Q_{x_1}(y + h+ k) = 0.\]
Since $Q_{x_1}$ is a homogeneous quadratic and $h,k$ (and $y$) were arbitrary, this last equation implies that $Q_{x_1}$ is in fact zero. 

Since no $x$ can be the $x_1$ term of more than $|\F|^{2n}$ additive quadruples, it follows that $Q_x = 0$ for at least $\frac{1}{100}\alpha^4 |\F|^n$ values of $x$.  On the other hand, the set of $x$ where $Q_x=0$ is a subspace of $\F^n$, and the claim is thus verified in this case.

\emph{Case 2.} At least half of the additive quadruples in $A$ are bad. Then (for example) there are at least $\frac{1}{100}\alpha^4 |\F|^{3n}$ quadruples $(x_1,x_2,x_3,x_4) \in A^4$ with $V_{x_1} \cap (V_{x_2} + V_{x_3}) \neq \{0\}$. Since the first three terms $x_1,x_2$ and $x_3$ of an additive quadruple determine the fourth, it follows easily that there is some choice of $x_2,x_3$ such that $V_{x_1} \cap (V_{x_2} + V_{x_3}) \neq \{0\}$ for at least $\frac{1}{100}\alpha^2|\F|^n$ values of $x_1$.  Since $V_{x_2} + V_{x_3}$ is a vector space of dimension at most $2r$, the claim follows in this case with $\alpha_0(\alpha,r) = \frac{1}{100}\alpha^2|\F|^{-2r}$.

We have verified the claim (with $\alpha_0(\alpha,r) = \frac{1}{100}\alpha^4 |\F|^{-2r}$, say) in all cases and hence the proposition is proved.\endproof

\emph{Remark.} An inspection of the argument reveals that the function $\alpha'(\alpha,r)$ in this proposition can be taken to have the form $(\alpha/C)^{C^r}$.\vspace{11pt}

Let us return now to \eqref{xr}, which stated that \[
 |\E_{y \in \F^N} e_{\F}((\phi(x)-\psi(x)) \cdot y -Q_x(y)-\Psi(0,y))| \gg_\delta 1 \]
for $\gg_\delta |\F|^n$ values of $x \in S$. In the subsequent discussion we passed to a further subset of $\gg_{\delta} |\F|^n$ values of $x$ for which $\rk(Q_x - Q_0) \ll_{\delta} 1$. Writing $Q'_x$ for the homogeneous quadratic part of $Q_x - Q_0$, we may use Proposition \ref{rank-line} to assert that there is some subspace $V \leq \F^N$, $\dim V \ll_{\delta} 1$, such that $Q'_x(y)$ is a quadratic function of the inner products $\langle v,y\rangle$, $v \in V$. The coefficients of this quadratic function vary \emph{linearly} in $x$, but this is unimportant.

By foliating into cosets of $V^{\perp}$, we may find a $1$-bounded function $F$ supported on some coset $t + V^{\perp}$ and a quadratic polynomial $\tilde\psi : \F^n \rightarrow \F^N$  such that
\[ \E_{y \in \F^N}F(y) e_{\F}((\phi(x) - \tilde \psi(x))\cdot y) \gg_{\delta} 1\] for $\gg_{\delta} |\F|^n$ values of $x \in S$. Note that the quadratic $\tilde \psi$ has been adjusted to take account for the possibility that $Q_x$ contains linear terms in $y$ (which also depend affine-linearly on $x$).

To conclude the argument we simply apply the Plancherel argument from Section \ref{ffc}. This tells us that there are $\ll_{\delta} 1$ values of $r$ for which\[ \E_{y \in \F^N} F(y) e_{\F}(r \cdot y) \gg_{\delta} 1.\] It follows from the pigeonhole principle that there is some $r$ such that $\phi(x) - \tilde\psi(x) = r$ for $\gg_{\delta} |\F|^n$ values of $x \in S$, which implies Proposition \ref{approx-hom-f5}.\endproof

\emph{Remark.} Because of the use of the rank reduction argument in the proof of Proposition \ref{rank-line}, the proof above does not seem to imply any implication between a conjectural polynomial version of Theorem \ref{conj-4}, and a polynomial version of Proposition \ref{approx-hom-f5}.  Also, we do not know if the implication can be reversed; the proof of Theorem \ref{conj-4} in \cite{btz,tz}, is somewhat different from the arguments in \cite{gowers:ap4,gowers-longaps,green-tao-u3inverse,samorodnitsky}, relying instead on ergodic theory and cohomological tools.

\appendix

\section{Deduction of PGI(3) from PFR}\label{gowi}

In this appendix we sketch how the polynomial Fre\u{\i}man-Ruzsa conjectures (Conjectures \ref{PFR-2}, \ref{PFR-1}) imply their respective polynomial inverse conjectures for the Gowers norms (Conjectures \ref{PGI-1}, \ref{PGI-2}).  Roughly speaking, the idea is to run the arguments in \cite{samorodnitsky} or \cite{green-tao-u3inverse} verbatim, but substituting the polynomial Fre\u{\i}man-Ruzsa conjectures in one key step of the argument where the usual inverse sumset theorems (basically, Theorem \ref{inv1} or \ref{inv2} respectively) are currently used instead.  It should be noted that the bulk of this implication is due to Gowers \cite{gowers:ap4,gowers-longaps}.

Our sketch will be somewhat brief and in particular we will assume familiarity with either \cite{samorodnitsky} or \cite{green-tao-u3inverse} as appropriate.  In the finite field case (i.e. the deduction of Conjecture \ref{PGI-1} from Conjecture \ref{PFR-2}) the modification is particularly straightforward; one simply repeats the argument in \cite{samorodnitsky}, but replacing \cite[Theorem 6.9]{samorodnitsky} (which is essentially Theorem \ref{inv1}) by Conjecture \ref{PFR-2} instead.  To spell out the steps in a little more detail, suppose that $K \geq 2$, and let $f: \F_2^n \to \C$ be a $K$-approximate quadratic: that is to say $\Vert f \Vert_{U^3(\F_2^n)} \geq 1/K$. By repeating the arguments up to and including \cite[Lemma 6.7]{samorodnitsky}, one can find a function $\phi: \F_2^n \to \F_2^n$ such that the set
$$
\{ (x,y) \in \F_2^n \times \F_2^n: \phi(x+y)=\phi(x)+\phi(y); |\widehat{f_x}(x)|, |\widehat{f_y}(y)|, |\widehat{f_{x+y}}(x+y)| \geq K^{-C} \}$$
has density $\geq K^{-C}$ in $\F_2^n \times \F_2^n$, where $f_x(y) := f(x+y) \overline{f(x)}$ and $\hat f(x) = \E_{y \in \F_2^n} f(y) (-1)^{x \cdot y}$ is the usual Fourier transform.  Now let
$$ A := \{x \in \F_2^n: |\widehat{f_x}(x)| \geq K^{-C} \}.$$
Arguing as in \cite[Section 6]{samorodnitsky}, but using Conjecture \ref{PFR-2} instead of \cite[Theorem 6.9]{samorodnitsky}, one finds a linear transformation $D: \F_2^n \to \F_2^n$ and $z \in \F_2^n$ such that $\phi(x) = D x + z$ for a proportion at least $c K^{-C}$ of all $x \in A$, and thus
$$ \E_{x \in \F_2^n} |\widehat{f_x}(D x + z)|^2 \geq K^{-C}.$$
By modulating $f$ by a suitable linear phase we may normalise so that $z=0$.  Continuing the argument in \cite[Section 6]{samorodnitsky} one concludes that the subspace $U := \{ x \in \F_2^n: D x = D^t x \}$ of $\F_2^n$ has density $\geq K^{-C}$, and so by further continuation of the argument one can find a symmetric transformation $B: \F_2^n \to \F_2^n$ with zero diagonal coefficients such that
$$ \E_{x \in \F_2^n} |\widehat{f_x}(Bx)|^2 \geq K^{-C}.$$
From the structure of $B$ one can $B = M + M^t$ for some transformation $M: \F_2^n \to \F_2^n$.  A little Fourier analysis then shows that the function $(-1)^{x \cdot Mx} f(x)$ has a $U^2(\F_2^n)$ norm of at least $K^{-C}$, and so has an inner product of at least $K^{-C}$ with a linear character, and Conjecture \ref{PGI-1} follows.

We turn now to the integer case, i.e. the deduction of Conjecture \ref{PGI-2} from Conjecture \ref{PFR-1}.  This requires a little more modification, because the arguments in \cite{green-tao-u3inverse} proceeded not via inverse sumset theorems, but instead via the (closely related) device of Bogulybov-type theorems\footnote{It is possible that polynomial variants of these Bogolyubov-type theorems also hold, but so far as we know conjectures of this type are strictly stronger than Conjectures \ref{PFR-2} and \ref{PFR-1}.}. We think, in particular of \cite[Lemma 6.3]{green-tao-u3inverse}).  However, as noted in \cite{gowers:ap4}, one could substitute inverse sumset theorems for Bogulybov-type theorems at this stage.

We turn to the details.  Let $K \geq 2$ and suppose that $f: \Z/N\Z \to \C$ is a $K$-approximate quadratic where, for simplicity, $N$ is odd (this in fact implies the general case, an exercise we leave to the reader).  Applying \cite[Proposition 5.4]{green-tao-u3inverse}, there is a set $H' \subset \Z/N\Z$ of size $|H'| \geq K^{-C} N$ and a function $\xi: H' \to \Z/N\Z$ whose graph
$$ \Gamma' := \{ (h,\xi_h): h \in H' \}$$
is such that $|9\Gamma'-8\Gamma'| \leq K^C N$, and such that
$$ |\hat f_h(\xi_h)| \geq K^{-C}$$
for all $h \in H''$, where $f_h(x) := f(x+h) \overline{f(x)}$ as before, and $\hat f(\xi) := \E_{\xi \in \Z/N\Z} f(x) e_N(x \xi)$ is the usual Fourier transform.
	
Applying \cite[Proposition 9.1]{green-tao-u3inverse}, one obtains a regular Bohr set $B_1 := B(S,\rho)$ with $|S| \leq K^C$, $\frac{1}{16} \leq \rho \leq \frac{1}{8}$ and $x_0, \xi \in \Z/N\Z$, as well as a locally linear function $M: B(S,\frac{1}{4}) \to \Z/N\Z$ such that
\begin{equation}\label{star-66}\E_{h \in B_1} 1_{H'}(x_0+h) 1_{\xi_{x_0+h} = 2Mh + \xi_0} \gg K^{-C}.\end{equation}
This was eventually used in \cite{green-tao-u3inverse} to deduce Theorem \ref{invu3}.  An inspection of that deduction reveals that the argument would also work just as well if the Bohr set $B(S,\rho)$ were replaced with a symmetric progression of dimension at most $K^C$ and cardinality at least $\exp(-K^C) N$.  Furthermore, if one could instead replace $B(S,\rho)$ with a progression of dimension at most $K^{o(1)}$ and cardinality at least $\exp(K^{-o(1)}) N$ then one could conclude Conjecture \ref{PGI-2} instead of Theorem \ref{invu3}.  Thus, our only task is to alter the argument of \cite[Proposition 9.1]{green-tao-u3inverse}, using the additional input of Conjecture \ref{PFR-2}, to obtain such a progression in place of $B(S,\rho)$.

By Conjecture \ref{PFR-1} $\Gamma'$ has large intersection with a translate of a symmetric generalised arithmetic progression $P$ of dimension at most $K^{o(1)}$ and cardinality at most $e^{K^{o(1)}} N$.  By \cite[Theorem 3.40]{tao-vu}, $P$ contains a \emph{proper} symmetric generalised arithmetic progression $P'$ 
$$P' = \{l_1 x_1 + \dots + l_d x_d : |l_i| \leq L_i\}$$  
in $\Z/N\Z \times \Z/N\Z$ of dimension $d \leq K^{o(1)}$ and volume at least $e^{K^{o(1)}} N$.  
The progression $P'-P'$ need not be a graph.  However, since $P'-P'+\Gamma' \subset 2P-P$ has size at most $e^{K^{o(1)}} N$, and $\Gamma'$ is a graph, we see that the intersection of $P'$ with the vertical axis $\{0\} \times \Z/N\Z$ has cardinality at most $e^{K^{o(1)}}$, thus $P'$ is in some sense ``almost a graph'' up to factors of $e^{K^{o(1)}}$.  Applying \cite[Lemma 8.3]{green-tao-u3inverse} one can then find a Bohr set $B(S,\frac{1}{4})$ in $\Z/N\Z$ with $|S| \leq K^{o(1)}$  such that $P-P \cap (\{0\} \times B(S,\frac{1}{4})) = \{0\}$.  In particular, the set $P'' := P' \cap (\Z/N\Z \times B(S,\frac{1}{8}))$ is a graph.  

One can write $P'' := \phi(B)$, where $B$ is the box $\{ (l_1,\ldots,l_d) \in \Z: |l_i| \leq L_i\}$ in $\Z^d$ and $\phi: \Z^d \to (\Z/N\Z \times \Z/N\Z)$ is the homomorphism $\phi(l_1,\ldots,l_d) := l_1 x_1 + \ldots + l_d x_d$.  Observe that $P'' = \phi(B \cap B(S',\frac{1}{8}))$ for some Bohr set $B(S',\frac{1}{8})$ in $\Z^d$. Applying Lemma \ref{prog-lem} $|S'|$ times we see that $B \cap B(S',\frac{1}{8})$ contains a symmetric generalised arithmetic progression $Q$ of dimension at most $K^{o(1)}$ and volume at least $e^{-K^{o(1)}} N$.  By shrinking $Q$ slightly we may in fact assume that $Q-Q \subset B \cap B(S',\frac{1}{8})$.  Then $\phi(Q-Q)$ is a graph, or equivalently that $\phi(Q)$ is Fre\u{\i}man isomorphic to its projection $\pi(\phi(Q))$ to the first factor $\Z/N\Z$ of $\Z/N\Z \times \Z/N\Z$.  Since $P'$ was proper, we see that $\pi(\phi(Q))$ is also proper.  We then conclude that
$$ \phi(Q) = \{ (x,Mx +\xi): x \in \pi(\phi(Q)) \}$$
where $\xi \in \Z/N\Z$, and $M: \pi(\phi(Q)) \to \Z/N\Z$ is locally linear.

As $Q$ is a progression, we can find $Q'-Q'$ inside $Q$ where $Q' \subset Q$ is another progression with dimension at most $K^{o(1)}$ and cardinality at least $e^{-K^{o(1)}} N$.  The set $\phi(Q')$ has relative density at least $e^{-K^{o(1)}}$ inside $P$, which has a doubling constant of at most $e^{K^{o(1)}}$, so by standard covering lemma arguments (see e.g. \cite[Lemma 2.14]{tao-vu}) one can cover $P$ by at most $e^{K^{o(1)}}$ translates of $\phi(Q')-\phi(Q') \subset \phi(Q)$.  In particular, by the pigeonhole principle, some translate of $\phi(Q)$ intersects $\Gamma'$ in at least $e^{-K^{o(1)}} N$ points.  If one then repeats the arguments used to prove \cite[Proposition 9.1]{green-tao-u3inverse} one obtains what was claimed, namely an analogue of \eqref{star-66} with $B(S,\rho)$ replaced by a progression of dimension $K^{o(1)}$ and size at least $\exp(-K^{o(1)})N$.

\section{Bohr sets in generalised progressions}\label{prog-app}

The aim of this appendix is to prove Lemma \ref{prog-lem}, the statement of which was as follows. 

\begin{lemma}\label{prog-lem} Let $\eps \in (0,1/2)$ be a real number.
Suppose that $P$ is a $d$-dimensional proper progression with sidelengths $N_1,\dots,N_d > C/\eps$ and that $\eta : P \rightarrow \R/\Z$ is a Fre\u{\i}man homomorphism which vanishes at some point of $P$. Then the set $\{x \in P : \Vert \eta(x)\Vert_{\R/\Z} \leq \eps\}$ contains a progression of dimension at most $d+1$ and size at least $(Cd)^{-d}\eps^{d+1}|P|$.
\end{lemma}
\proof The progression $P$ is an affine image of some box $[N_1] \times \dots \times [N_d]$, and the lift of $\eta$ to this box is an affine map of the form $x \rightarrow \alpha_1 x_1 + \dots + \alpha_d x_d + \beta$. Henceforth we abuse notation by identifying $P$ with the box $[N_1] \times \dots \times [N_d]$. We are told that there is a point $x^*$ such that $\eta(x^*) = 0$.  By reparametrising $P$ if necessary, we may assume that $x^*$ is in the same quadrant of $P$ as the origin, thus $x^* \in [N_1/2] \times \dots \times [N_d/2]$. It turns out to be inconvenient later on if $x^*$ is too close to the boundary of $P$, so we begin with a preliminary argument to find a point $x^{**}$ which is deeper in the interior of $P$ than $x^*$, and at which $\eta$ is still small.  To do this consider some $m := \lceil 2/\eps \rceil + 1$ points $x_1,\dots,x_m \in P$ such that the $j$th coordinate of $x_i$ is roughly $i N_j/3m$. By the pigeonhole principle there must be some pair of indices $s < t$ such that $\Vert \eta(x_t - x_s)\Vert_{\R/\Z} \leq \eps/2$, and then the point $x^{**} := x^* + x_t - x_s$ will have the property that all of its coordinates lie between $\eps N_j/10$ and $(1 - \eps/10)N_j$ (note that we implicitly used here the fact that $N_j > C/\eps$).

Let us now recentre so that $x^{**}$ is at the origin. Since $x^{**}$ was chosen to be somewhat central to $P$, the progression $P$ certainly contains the symmetric progression $P' := \prod_{j=1}^d [-N'_j,N'_j]$ in this new coordinate system, where $N'_j := \eps N_j/10$. Henceforth we work entirely in this new coordinate system and with this new progression $P'$. The Fre\u{\i}man homomorphism $\eta: P' \to \R/\Z$ now takes the form $\eta(x) = \alpha_1 x_1 + \dots + \alpha_d x_d + \beta$ where $\Vert \beta \Vert_{\R/\Z} \leq \eps/2$, and we may of course assume that $0 \leq \alpha_j < 1$ for each $j$.

At this point, one could conclude the argument (with worse bounds than claimed) using \cite[Lemma 4.20, Lemma 4.22]{tao-vu}, because the set where $\eta$ is small is essentially a Bohr set in $P'$.  To get the sharper bounds claimed in the theorem, we use a well-known lemma of Ruzsa \cite[Theorem 3.1]{ruzsa-freiman}, in which the structure of Bohr sets was elucidated using the geometry of numbers.

\begin{lemma}\label{bohr-prog} Suppose that $M \geq 1$ is an integer, that $r_1,\dots,r_d$ are residues $\mdlem{M}$ such that $\hcf(r_1,\dots,r_k,M) = 1$, and that $\eps _1,\dots\eps_d \in (0,1/2)$ are real numbers. Then  the \emph{Bohr set} 
\[ B(r_1,\dots,r_d;\eps_1,\dots,\eps_d) := \{x \in \Z/M\Z : \Vert r_1x/M\Vert_{\R/\Z}  \leq \eps_1,\dots,\Vert r_d x/M\Vert_{\R/\Z} \leq \eps_d\}\] contains a $d$-dimensional progression \textup{(}that is, the image of a box under an affine map from $\Z^d$ to $\Z/M\Z$\textup{)} of cardinality at least $d^{-d}\eps_1\dots \eps_d M$.\endproof
\end{lemma}

Let $M_1 \geq \dots \geq M_d$ be a very large odd coprime integers and set $M := M_1\dots M_d$. Set $r_j := M_{j+1}\dots M_d$ for $j = 1,\dots,d-1$ and $r_d := 1$. For each $j= 1,\dots,d$ choose an integer $s_j$, $0 \leq s_j < M_j$, such that $|s_j/M_j - \alpha_j| \leq 1/M_j$. Set $r_{d+1} := r_1 s_1 + \dots + r_d s_d$. Finally, set $\eps_j := N'_j/2M_j$ for $j = 1,\dots,d$ and $\eps_{d+1} := \eps/4$. Our contention is that the Bohr set $B' = B(r_1,\dots,r_{d+1};\eps_1,\dots,\eps_{d+1})$ is contained in a set which is Fre\u{\i}man isomorphic to $\{x \in P: \Vert \eta(x) \Vert_{\R/\Z} \leq \eps\}$, at which point Lemma \ref{prog-lem} follows easily from Lemma \ref{bohr-prog}.
To begin with we show that the Bohr set $B = B(r_1,\dots,r_d;\eps_1,\dots,\eps_d)$ is contained in a Fre\u{\i}man-isomorphic copy of $P$. Suppose that $x \in \Z/M\Z$ lies in $B = B(r_1,\dots,r_d;\eps_1,\dots,\eps_d)$. If $x \in \Z/M\Z$, we may write
\[ x = x_1 + x_2 M_1 + \dots + x_d M_1 \dots M_{d-1}\] for unique integers $x_1,\dots,x_d$ with $|x_j|  < M_j/2$. Observe that 
\[ \frac{r_1 x}{M} \equiv \frac{x_1}{M_1} \md{1},\qquad \frac{r_2 x}{M} \equiv \frac{x_1}{M_1M_2} + \frac{x_2}{M_2} \md{1},\] and so on. If $x \in B$ then these may be applied in succession to obtain $\Vert r_1 x/M - x_1/M_1 \Vert_{\R/\Z} = 0$, 
\begin{equation}\label{b2}  \Vert \frac{r_2 x}{M} - \frac{x_2}{M_2} \Vert_{\R/\Z} \leq \frac{\eps_1}{M_2},\qquad  \Vert \frac{r_3 x}{M} - \frac{x_3}{M_3} \Vert_{\R/\Z} \leq \frac{\eps_1}{M_2M_3} + \frac{\eps_2}{M_3},\end{equation} and so on. 
If the $M_j$ are chosen appropriately (with $M_1$ much bigger than $M_2$ and so on) this implies that $\Vert x_j/M_j\Vert \leq 2\eps_j$ for $j = 1,\dots,d$, which implies that $|x_j| \leq N'_j$ for all $j$. 

Now we have
\begin{equation}\label{b3} \Vert \frac{s_1 x_1}{M_1} + \dots + \frac{s_dx_d}{M_d} - \eta(x)\Vert_{\R/\Z}  \leq  |\frac{s_1}{M_1} - \alpha_1| |x_1| + \dots + |\frac{s_d}{M_d} - \alpha_d| |x_d|  \leq \frac{N'_1}{M_1} + \dots + \frac{N'_d}{M_d} \leq \frac{\eps}{8},
\end{equation}
provided that the $M_j$ are chosen large enough in terms of $N'_1,\dots,N'_d$ and $\eps$.

Furthermore the inequalities \eqref{b2} imply that $\Vert s_1r_1 x/M - s_1x_1/M_1 \Vert = 0$,
\[ \Vert \frac{s_2r_2 x}{M} - \frac{s_2x_2}{M_2} \Vert_{\R/\Z} \leq \frac{s_2\eps_1}{M_2} \leq \eps_1,\] \[ \Vert \frac{s_3r_3 x}{M} - \frac{s_3x_3}{M_3} \Vert_{\R/\Z} \leq \frac{s_3\eps_1}{M_2M_3} + \frac{s_3\eps_2}{M_3} \leq \frac{\eps_1}{M_2} + \eps_2,\] 
and so on.
Adding, we clearly obtain
\[ \Vert \frac{r_{d+1} x}{M} - \frac{s_1 x_1}{M_1} - \dots - \frac{s_d x_d}{M_d} \Vert_{\R/\Z} \leq \frac{\eps}{8}\] provided that the $M_i$ are selected to be large enough. 

Combining this with \eqref{b3}, we obtain 
\[ \Vert \frac{r_{d+1} x}{M} - \alpha_1 x_1 - \dots - \alpha_d x_d\Vert_{\R/\Z} \leq \frac{\eps}{4},\] and hence if $x \in B$ we certainly have $\Vert \alpha_1 x_1 + \dots + \alpha_d x_d\Vert_{\R/\Z} \leq \eps/2$ and hence $\Vert \eta(x) \Vert_{\R/\Z} \leq \eps$.\endproof

\providecommand{\bysame}{\leavevmode\hbox to3em{\hrulefill}\thinspace}
\providecommand{\MR}{\relax\ifhmode\unskip\space\fi MR }
\providecommand{\MRhref}[2]{%
  \href{http://www.ams.org/mathscinet-getitem?mr=#1}{#2}
}
\providecommand{\href}[2]{#2}

     \end{document}